\documentclass[10pt]{article}
\usepackage{amssymb,amsmath, amsfonts, amsthm, amscd}

\title{Notes on motives in finite characteristic}

\author {Maxim Kontsevich}

\begin{document}

\maketitle

\newtheorem{dfn}{Definition}
\newtheorem{thm}{Theorem}

\newtheorem{lmm}{Lemma}
\newtheorem{rmk}{Remark}
\newtheorem{prp}{Proposition}
\newtheorem{conj}{Conjecture}
\newtheorem{exa}{Example}
\newtheorem{cor}{Corollary}
\newtheorem{que}{Question}

\newcommand{\G}{\Gamma}
\newcommand{\A}{{\mathbb A}}
\newcommand{\Z}{{\mathbb Z}}
\newcommand{\C}{{\mathbb C}}
\newcommand{\Q}{{\mathbb Q}}
\newcommand{\R}{{\mathbb R}}

\newcommand{\GG}{{\mathbb G}}
\renewcommand{\P}{{\mathbb P}}
\newcommand{\Qb}{{\overline{\Q}}}
\newcommand{\F}{{\mathbb F}}
\newcommand{\CC}{{\mathcal C}}
\renewcommand{\AA}{{\mathcal A}}
\renewcommand{\k}{{\bf  k}}
\newcommand{\kb}{{\overline{\k}}}

\newcommand{\E}{{\mathcal E}}

\newcommand{\B}{{\mathcal B}}

\newcommand{\Fb}{{\overline{\F}}}
\newcommand{\FF}{{\op{Fun}^{poor}}}
\newcommand{\1}{{\bf 1}}

\newcommand{\spec}{{\mathtt {Spec}}}
\newcommand{\Gal}{{\op{Gal}}}

\newcommand{\op}[1]{\operatorname{#1}}

\newcommand{\epi}{\twoheadrightarrow}
\newcommand{\mono}{\hookrightarrow}

\newcommand{\tto}{\longrightarrow}
\newcommand{\ltto}{\longleftarrow}
\newcommand\uhom{{\underline{Hom}}}
\newcommand{\umap}{{\underline{Map}}}
\renewcommand\O{{\cal O}}

\newcommand{\Fr}{\op{Fr}}

\ \ \ \ \ \ {\it To Yuri Manin on the occasion of 70-th birthday, with admiration.}

\section*{Introduction and an example}

These notes grew from an attempt to interpret a formula of Drinfeld (see \cite{Dr}) enumerating the absolutely irreducible local systems of rank 2 on algebraic curves over 
finite fields, obtained as
 a  corollary of the Langlands correspondence for $GL(2)$ in the functional field case, and of the trace formula.
 
Let $C$ be a smooth projective geometrically connected
 curve defined over a finite field $\F_q$, with a base point $v\in C(\F_q)$. The geometric
 fundamental group $\pi_1^{geom}(C,v):=\pi_1(C\times_{\spec \,\F_q}\spec\,\Fb_q,v)$ is a profinite group
on which the Galois group $\widehat{\Z}=\Gal(\Fb_q/\F_q)$ (with the canonical generator $\Fr:=\Fr_q$)
 acts. In what follows we will omit the base point from the notation.

\begin{thm}{\bf (Drinfeld)} Under the above assumptions, 
for any integer $n\ge 1$ and any prime  $l\ne \mbox{char}(\F_q)$ 
 the  set of fixed points
$$X_n^{(l)}:=\left(\op{IrrRep}\,(\pi_1(C\times_{\spec \,\F_q}\spec\,\Fb_q)\to
GL(2,\Qb_l))/\mbox{conjugation}\right)^{\Fr^n}$$
 is finite. 
Here $\op{IrrRep}\,(\dots)$ denotes the set of conjugacy classes of irreducible continuous
$2$-dimensional representations of $\pi_1^{geom}(C)$ defined over  finite extensions of $\Q_l$.
Moreover, there exists a finite collection $(\lambda_i)\in  \Qb^\times$
 of algebraic integers, and signs $(\epsilon_i)\in\{-1,+1\}$ depending only on $C$,
such that for any $n,l$ one has an equality
 $$ \# X_n^{(l)}=\sum_i \epsilon_i \lambda_i^n
\,\,.$$
\end{thm}

From the explicit formula which one can extract from \cite{Dr}
 one can see that numbers $\lambda_i$ are $q$-Weil algebraic integers whose norm for any  embedding $\Qb\mono \C$ belongs to
$q^{\frac{1}{2}\Z_{\ge 0}}$.
 Therefore, the number of elements of $X_n^{(l)},\,n=1,2,\dots$  looks like the number of $\F_{q^n}$-points
 on some variety over $\F_q$. The largest exponent is $q^{4g-3}$, which indicates that this variety
 has dimension $4g-3$. A natural guess is that it is closely related to the moduli
space of stable bundles of rank $2$ over $C$. At least the
dimensions coincide, and Weil numbers which appear
 are essentially the same, they are 
products of the eigenvalues of Frobenius acting on the motive defined by the first cohomology of $C$.

 The Langlands correspondence identifies $X_n^{(l)}$ with the set of $\Qb_l$-valued unramified cuspidal automorphic forms 
 for the adelic group $GL(2,\A_{\F_q(C)})$. These forms are eigenvectors of a collection of commuting matrices (Hecke operators)
  with integer coefficients.  Therefore, for a given $n\ge 1$ one can identify\footnote{It is expected that all 
  representations from $X_n^{(l)}$ 
  are  motivic, i.e. they arise from projectors with coefficients in $\Qb$ acting
  on $l$-adic cohomology of certain projective varieties defined over the field of rational functions $\F_{q^n}(C)$.}
    all sets $X_n^{(l)}$ for various primes $l$ with one set $X_n$ endowed with an action
   of the absolute Galois group $\Gal(\Qb/\Q)$, extending the obvious actions of $\Gal(\Qb_l/\Q_l)$ on $X_n^{(l)}$.

These days the Langlands correspondence in the functional field case is established for all the groups $GL(N)$
 by L.~Lafforgue. To my knowledge, almost no attempts were made to extend Drinfeld's calculation
 to the case of higher rank, or even to the $GL(2)$ case
with non-trivial ramification.

It is convenient
 to take the inductive limit $X_\infty:=\varinjlim X_n,\,\,\,X_{n_1}\hookrightarrow X_{n_1 n_2}$ which is an
infinite countable set endowed with an action of the
product\footnote{One can replace
 $\Gal(\Qb/\Q)$  by its quotient $\Gal(\Q^{q-Weil}/\Q)$ where $\Q^{q-Weil}\subset\Qb$
is CM-field generated by all $q$-Weil numbers.}
 $$\Gal(\Qb/\Q)\times \Gal(\Fb_q/\F_q)\,\,\,.$$
The individual set $X_n$ can be reconstructed from this datum as the set of fixed points of  $\Fr^n\in \Gal(\Fb_q/\F_q)$.

In spite of the numerical evidence, it would be too naive to expect  a natural identification
of $X_\infty$  with the set of $\Fb_q$-points of an algebraic variety defined over $\F_q$, as
there is no obvious mechanism producing a non-trivial $\Gal(\Qb/\Q)$-action on the latter.
 
 The main question addressed here is 
\begin{que} Does there exist some alternative way to construct the set $X_\infty$
 with the commuting action of two Galois groups?
\end{que}

In the present notes I will offer three different hypothetical constructions.
The first construction comes from the analogy between the Frobenius acting on $\pi_1^{geom}(C)$ and an element of the mapping class
 group acting on the fundamental group of a closed oriented surface,
 the second one
 is almost tautological and arises from the contemplation on the shape of explicit
 formulas for Hecke operators (see an example in Section 0.1),
the third one is based on an analogy with lattice models in statistical physics.

I propose several conjectures, which should be better considered
as  guesses 
 in the first and in the third part, as there is almost no experimental
evidence in their favor.
In a sense, the first and the third part should be regarded as
  science fiction, but even if the appropriate conjectures are wrong (as I strongly suspect),
 there should be some grains of truth in them.

On the contrary,  I feel quite confident that the conjectures made
in the second part  are essentially true,
 the output is a higher-dimensional generalization of the Langlands correspondence in the functional field case.
  At the end of the second part I will show how to make a step in the arithmetic direction, extending
   the formulas to the case of an arbitrary local field.

In the fourth part  I will describe briefly a similarity between a  
modification of the category of motives
based on non-commutative geometry, and two other categories introduced in the second and the third part.
Also I will make a link between the proposal based on polynomial dynamics and the one based
on lattice models.
 
 Finally, I apologize to the reader that the formulas in Sections 0.1 and 1.3 are given without explanations,
 this is the result of my poor knowledge of the representation theory. The formulas were polished with the help
 of computer.

{\it Acknowledgements:} I am grateful to many people for useful discussions, especially to Vladimir Drinfeld
 (on the second part of this paper),
to  Laurent Lafforgue who told me about \cite{Dr} and explained some basic stuff about automorphic forms,
 to  Vincent Lafforgue who proposed an argument for the Conjecture 5 based on lifting,  to James Milne
 for consultations about the category of motives over a finite field,
to Misha Gromov who suggested the idea to
 imitate Dwork's methods in the third part, also to
  Mitya Lebedev, Sasha Goncharov, Dima Grigoriev,
 Dima Kazhdan, Yan Soibelman and Don Zagier. 
 Also I am grateful to the referee for useful remarks and corrections.

\subsection{An explicit example}

Here we will show explicit formulas for the tower $(X_n)_{n\ge 1}$ in the simplest truly non-trivial case.
Consider the affine curve $C=\P^1_{\F_q}\setminus\{0,1,t,\infty\}$
for a given element $t\in \F_q\setminus\{0,1\}$. We are
interested in motivic local
 $SL(2)$-systems on $C$ with tame non-trivial unipotent monodromies around all punctures $\{0,1,t,\infty\}$.
 
 A lengthy calculation lead to the following explicit formulas\footnote{I was informed by V.~Drinfeld that a similar calculation for the case of $SL(2)$ local systems on
$\P^1_{\F_q}\setminus \{4\mbox{ points}\}$, with tame non-trivial
{\it{semisimple}} monodromy around punctures was performed few years ago by Teruji Thomas.} for the Hecke
 operators for cuspidal representations. In what follows we assume $\op{char} \F_q\ne 2$. The Hecke operators act on 
 the spaces of functions on certain double coset space for the adelic group, which can be identified
 with the set of equivalence classes of vector bundles of rank 2 over $\P^1_{\F_q}$ 
 together with a choice of one-dimensional subspaces of fibers at $\{0,1,t,\infty\}$.
 This double coset space is infinite, but the eigenfunctions of Hecke operators corresponding to cuspidal representations have finite support which one can bound a priori.

For any $x\in \F_q$ the Hecke operator  $T_x$ (on cuspidal forms) can be written as an integral $q\times q$
 matrix whose rows and columns are labelled by elements of $\F_q$,
 (i.e. $T_x\in \op{Mat}(\F_q\times \F_q;\Z)$). Coefficients of $T_x$ are given by the formula
$$(T_x)_{yz}:=2-\#\{w\in \F_q|\,w^2=f_t(x,y,z)\}+\mbox{(correction term)}$$
where $f_t(x,y,z)$ is the following universal polynomial with integral coefficients:
$$f_t(x,y,z):=(xy+yz+zx-t)^2+4xyz(1+t-(x+y+z))\,\,\,.$$
The correction term is equal to
$$-\left\{
{\begin{aligned} q&+1 && x=y\in\{0,1,t\}\\
                     &1        && x=y\notin \{0,1,t\}\\
                     &0        && x\ne y
\end{aligned}}
 \right.
+\left\{ \begin{aligned}
 &q &&\mbox{ if } x\notin\{0,1,t\}\,\,\mbox{ and }
\left\{\begin{aligned} y&=\frac{t}{x}, && z=0\\
                      y&=\frac{t-x}{1-x}, &&z=1\\
                      y&=\frac{t(1-x)}{t-x}, &&z=t
\end{aligned}\right.
 \\
       & 0 &&\mbox{ otherwise}
          \end{aligned}
\right.$$

Operators $(T_x)_{x\in \F_q}$ satisfy the following properties:
\begin{enumerate}
\item $[T_{x_1},T_{x_2}]=0$,
\item $\sum_{x\in
\F_q}T_x=\1=id_{\Z^{\F_q}} $,
 \item $T_x^2=\1$ for
$x\in\{0,1,t\}$, moreover $\{\1,T_0,T_1,T_t\}$ form a group under the multiplication, isomorphic to $\Z/2\Z\oplus
\Z/2\Z$,
\item for any $x\notin \{0,1,t\}$ the spectrum of $T_x$ is real and belongs to $[-2\sqrt{q},+2\sqrt{q}]$,
 any element of $\spec(T_x)$ can be written as $\lambda+\overline{\lambda}$ where $|\lambda|=\sqrt{q}$
is a $q$-Weil number,
\item for any $\xi=\lambda+\overline{\lambda}\in \spec(T_x)$ and any integer $n\ge 1$
the spectrum of the matrix $T_x^{(n)}$ corresponding to $x\in \F_q\subset\F_{q^n}$ (if we pass
 to the extension $\F_{q^n}\supset \F_q$) contains the 
element $\xi^{(n)}:=\lambda^n+\overline{\lambda}^n$,
\item the vector space generated by $\{T_x\}_{x\in \F_q}$ is closed
under the product, the multiplication table is
$$T_x\cdot T_y=\sum_{z\in \F_q} c_{xyz}T_z\,\mbox{ where }c_{xyz}= (T_x)_{yz} \,\,\,.$$
\end{enumerate}

Typically (for ``generic'' $t,x$) the characteristic polynomial of $T_x$ splits into the product of
$4$ irreducible polynomials of almost the same degree. The splitting is not surprising,
 as we have a group\footnote{This is the group of automorphisms of $\P^1\setminus \{4\mbox{ points}\}$
 for
the generic cross-ratio.}
 of order 4 commuting with all operators $T_x$ (see property 3).
 Computer experiments indicate that the Galois groups of these polynomials (considered as permutation groups)
tend to be rather large,
 typically the full symmetric groups if  $q$ is prime, and the corresponding number fields have huge
factors in the prime decomposition of the discriminant.

Notice that in the theory of automorphic forms one usually
deals with infinitely many commuting Hecke
 operators corresponding to all places of the global field, i.e. to closed points of $C$ (in other
 words, to orbits of $\Gal(\Fb_q/\Fb_q)$ acting on $C(\Fb_q)$).
Here we are writing formulas only for the points defined over $\F_q$. The advantage of our example is that the number
 these operators coincides with the size of Hecke matrices, hence one can try to write formulas for structure constants, which by luck turn out to coincide
  with the matrix coefficients of matrices $T_x$ (property 6).

\section{First proposal: algebraic dynamics}

As was mentioned before, it is hard to imagine a mechanism for a non-trivial action
of   the absolute Galois group of $\Q$ on the set of points of a variety over a finite field. One can try to exchange the roles of fields $\Q$ and $\F_q$.
The first proposal is the following one:

\begin{conj} For a tower $(X_n)_{n\ge 1}$ arising from automorphic forms (or from motivic local systems on
curves), as defined in the Introduction,
 there exists a variety $X$ defined over $\Q$
and a map $F:X\to X$ such that there is a family
 of  bijections 
$$X_n \simeq (X(\Qb))^{F^n},\,\,\,n\ge 1$$
covariant with respect to $\Gal(\Qb/\Q)\times \Z/n\Z$ actions, and with respect to inclusions $X_{n_1}\subset X_{n_1 n_2}$
for integers $n_1,n_2\ge 1$.
\end{conj}

\subsection{The case of $GL(1)$}

Geometric class field theory gives a description of the sets $(X_n)_{n\ge 1}$ 
in terms of the Jacobian  of $C$:
$$ X_n=(\op{Jac}_C(\F_{q^n}) )^\vee(\Qb)
= \op{Hom}(\op{Jac}_C(\F_{q^n}),\Qb^\times)\,\,.  $$

The number of elements of this set is equal to 
$$\# \op{Jac}_C(\F_{q^n})\,=\det(\Fr_{H^1(C)}^n-\1)$$ where $\Fr_{H^1(C)}$ is
the Frobenius operator acting on, say, $l$-adic first cohomology group of $C$.

One can propose a blatantly non-canonical candidate for the corresponding dynamical system $(X,F)$.
Namely, let us choose a semisimple $(2g\times 2g)$ matrix $A=(A_{i,j})_{1\le i,j\le 2g}$ (where $g$ is the
genus of $C$)
with coefficients
 in $\Z$, whose characteristic polynomial is equal to the characteristic polynomial of $\Fr_{H^1(C)}$.
 Define $X/\Q$ to be the standard $2g$-dimensional torus $\GG_m^{2g}=\op{Hom}(\Z^{2g},\GG_m)$,
 and the map $F$ to be the dual to the map $A:\Z^{2g}\to\Z^{2g}$:
    $$F(z_1,\dots,z_{2g})=(\prod_i z_i^{A_{i,1}},\dots,
\prod_i z_i^{A_{i,2g}})\,\,\,.$$

Moreover, one can choose $A$ in such a way that 
$$F^*\omega= q\omega \,\,\,\,\,\,\mbox{ where }\omega=\sum_{i=1}^g \frac{dz_i}{z_i}\wedge
\frac{dz_{g+i}}{z_{g+i}}\,\,\,.$$

On the set of fixed points of $F^n$ act simultaneously $\Gal(\Qb/\Q)$ (via the cyclotomic quotient) and
 $\Z/n\Z$ (by powers of $F$). Nothing  contradicts to the existence of an equivariant isomorphism between
 two towers of finite sets.

\subsection{Moduli of local systems on surfaces}

 One can interpret the scheme $\GG_m^{2g}$ as the moduli space of rank $1$ local systems on a oriented 
closed topological surface $S$ of genus $g$, the form $\omega$ is the natural symplectic form on this moduli space.

In general, for any $N\ge 1$,
 one can make an analogy between the action of Frobenius $\Fr$ on the the set of $l$-adic irreducible 
 representations of $\pi_1^{geom}(C)$ of rank $N$, and the action
 of the isotopy class of a homeomorphism $\varphi:S\to S$ on the set of 
 irreducible complex representations of $\pi_1(S)$ of the same rank. Sets of representations
  are similar to each other, as it is  known
  that the maximal quotient of $\pi_1^{geom}(C)$ coprime to $q$ is isomorphic to the analogous quotient of the 
  profinite completion $\widehat{\pi}_1(S)$ of $\pi_1(S)$.
  Also, if we assume that there are only finitely many fixed points of $\varphi$ acting on 
  $$\op{IrrRep}(\pi_1(S)\to
GL(2,\C))/\mbox{conjugation}$$
then the sets
$$X^{(l)}:=\left(\op{IrrRep}(\widehat{\pi}_1(C\times_{\spec \,\F_q}\spec\,\Fb_q)\to
GL(2,\Qb_l))/\mbox{conjugation}\right)^\varphi$$
 do not depend on the prime $l$ for $l$ large enough.

All this leads to the following conjecture (which is formulated a bit sloppily), a strengthening of Conjecture 1:

\begin{conj} For any smooth compact geometrically connected curve $C/\F_q$ of genus $g\ge 2$ there exists an endomorphism $\Phi_C$ of the tensor category 
of finite-dimensional complex local systems on $S$ such that 
\begin{itemize}
\item $\Phi_C$ is {\it algebraic and defined over } $\Q$, in the sense that  it acts on the moduli stack of irreducible 
local systems of any given rank $N\ge 1$ 
 by a rational map defined over $\Q$,
 \item $\Phi_C$ multiplies the natural symplectic form on the moduli space of irreducible local systems of rank $N$ by the constant $q$,
 \item for every $n,N\ge 1$ there exists an identification of  the set of isomorphism classes of irreducible motivic local systems of rank $N$ on 
 $C\times_{\spec \F_q}\spec\,\Fb_q$
 invariant under $\Fr^n$, with the set of isomorphism classes of $\Qb$-local systems  of rank $N$ on $S$ invariant under $\Phi_C^n$, 
 compatible with
  the relevant Galois symmetries and tensor constructions.
  \end{itemize}
  \end{conj}
  
  One can not expect that $\Phi_C$ comes from an actual endomorphism $\varphi$ of the fundamental group $\pi_1(S)$, as it is known that for $g\ge 2$
  any such $\varphi$ is necessarily an automorphism. That is a rationale for replacing a putative endomorphism of $\pi_1(S)$ by a more esoteric 
  endomorphism of the tensor category of its finite-dimensional representations.

\subsubsection{Example: $SL(2)$-local systems on sphere with 3 punctures}

A generic $SL(2,\C)$-local system on $\C P^1\setminus \{0,1,\infty\}$ is uniquely determined by $3$ traces of monodromies around punctures.
  A similar statement holds for $l$-adic local systems with  tame monodromy in the case of finite characteristic.
   Motivic local systems correspond to the case when all the eigenvalues of the monodromies around punctures are roots of unity, i.e. when the traces
    of monodormies are twice cosines of rational angles.
     This leads to the following prediction:
     $$X=\A^3, \,\,\,F(x_1,x_2,x_3)=(T_q(x_1),T_q(x_2),T_q(x_3))$$
     where $T_q\in \Z[x]$ is the $q$-th Chebyshev polynomial,
     $$T_q(\lambda+\lambda^{-1})=\lambda^q+\lambda^{-q}\,\,\,.$$
     In this case the identifications between the fixed points of $F$ and motivic local systems on $\P^1_{\F_q}\setminus \{0,1,\infty\}$
      exist, and can be extracted form the construction of these  local systems (called hypergeometric)
      as summands in certain direct images of abelian local systems
       (analogous the classical integral formulas for hypergeometric functions). The identification is ambiguous, it depends
        on a choice of a group embedding $\Fb_q^\times \mono \Q/\Z $.

\subsection{Equivariant bundles and Ruelle-type zeta-functions}

The analogy with an element of the mapping class group acting on surface $S$ suggest the following addition to the Conjecture 1.
 Let us fix the curve $C/\F_q$ and the rank $N\ge 1$ of local systems under the consideration.
 For a given point $x\in C(\F_q)$ we have a sequence of  Hecke operators $T_x^{(n)}$ associated with curves 
 $C\times_{\spec \,\F_q}\spec\,\F_{q^n}$.
 The spectrum of $T_x^{(n)}$ is a $\Qb$-valued function on $X_n$, i.e. according to Conjecture 1,  a function on 
 the set of fixed points of $F^n$.
  We expect that the collection of these functions for $n=1,2\dots$  comes from a $F$-equivariant vector bundle   on $X$.
  \begin{conj} Using the notations of  Conjecture 1, for given $x\in C(\F_q)$ there exists a
  pair $(\E,g)$ 
 where $\E$ is a vector bundle on $X$ of rank $N$ together with an isomorphism $g:F^*\E\to \E$ (defined over $\Q$), such that
  the eigenvalue of $T^{(n)}_x$ at the point of its spectrum corresponding to $z\in X_n$ coincides with
  $$ \op{Trace}(\E_z= \E_{F^{n}(z)}\to \dots\to \E_{F(z)}\to\E_z)$$
  where arrows are isomorphisms of fibers of $\E$ coming from $g$.
  \end{conj}

  In particular, one can ask for an explicit formula for the $F$-equivariant bundle $\E$ in the case
   of $SL(2)$-local systems on the sphere with 3 punctures where we have an explicit candidate for $(X,F)$.

 In the limiting most simple non-abelian case when the monodromy is
unipotent around 2 punctures, and arbitrary semisimple around the
third puncture, one can make the above question completely
explicit:
\begin{que} For a given $x\in \F_q\setminus\{0,1\}$, does there exist a rational function $R=R_x$ on $\C
P^1$ with values in $q^{1/2}SL(2,\C)\subset\op{Mat}(2\times 2, \C)$ which has no
singularities on the set
$$ \left(\cup_{n\ge 1}\{z\in\C |
z^{q^n-1}=1\}\right)\setminus\{1\}$$ such that for any $n\ge 1$
two sets of complex numbers (with multiplicities):
$$ X_n:=\left\{ \sum_{y\in \F_{q^n}\setminus\{0,1,x\}}
\chi\left(\frac{y(1-xy)}{1-y}\right)\left|\,\,\chi:\F_{q^n}^\times
\to \C^\times ,\,\,\chi\ne 1\right.\right\}$$ where $\chi$ runs
through all non-trivial multiplicative characters of $\F_{q^n}$,
and
$$ X'_n:=\left\{ \op{Trace}\left(R(z)R(z^q)\dots R(z^{q^{n-1}})\right)|\,z^{q^n-1}=1,\,z\ne 1\right\}$$
coincide?
  \end{que}

Elements of the set $X_n$ are real numbers of the form
$\lambda+\overline{\lambda}$ where $\lambda\in\Qb$ is a $q$-Weil
number with $|\lambda|=q^{1/2}$. Therefore it is natural to expect
that $R(z)$ belongs to $q^{1/2} SU(2)$ if $|z|=1$.

The Galois symmetry does not forbid for the function $R$ (as a rational
function with values in $(2\times 2)$-matrices) to be defined over $\Q$, after the conjugation
by a constant matrix. Moreover,  the 
existence of such a function over $\Q$  leads  to certain choice of generators of
 the multiplicative groups $\left(\F_{q^n}^\times\right)$ for all
 $n\ge 1$, well-defined modulo the  action of Frobenius $\Fr_q$
  (the Galois group $\Gal(\F_{q^n}/\F_q)=\Z/n\Z$), as in a sense we identify
     roots of unity in $\C$ and   multiplicative characters of $\F_{q^n}$. In particular,
      there will be a {\it canonical} irreducible polynomial over $\F_q$ of degree $n$ for every $n\ge 1$. This is
  something almost too good to be true.

\subsubsection{Reminder: Trace formula and Ruelle-type zeta-function}

Let $X$ be now a smooth {\it proper} variety (say, over $\C$), endowed with a map $F:X\to X$, and
$\E$ be a vector bundle
on $X$ together with a morphism (not necessarily invertible) $g:F^* \E\to \E$.
 Let us assume that for any $n\ge 1$ all fixed points $z$ of $F^{n}:X\to X$ are isolated and {\it non-degenerate},
  i.e.  the tangent map
       $$ (F^{n})'_z: T_z X\to T_z X$$
        has no nonzero invariant vectors (in other words,  all eigenvalues of $ (F^{n})'_z$ are not equal to $1$).
    Then one has the following identity (Atiyah-Bott fixed point formula):
    $$\sum_{v\in X:F^{n}(z)=z} \frac
    {\op{Trace}(\E_z=\E_{F^{n}(z)}\to  \dots\to\E_z)} {\det (1-(F^{n})'_z)}=$$
    $$=
    \op{Trace}((g_*\circ F^*)^n:H^\bullet(X,\E)\to H^\bullet(X,\E))$$
    The trace in the r.h.s. is understood in the super sense, as the alternating sum of the ordinary traces
    in individual cohomology spaces.

    If one wants to eliminate the determinant factor in the denominator in the l.h.s., one should replace
     $\E$ by the superbundle $\E\otimes \wedge^\bullet (T^*_X)$.

     The trace formula implies that the series in $t$
     $$\exp\left(
     -\sum_{n\ge 1} \frac{t^n}{n}\sum_{z\in X:F^{n}(z)=z} \frac
    {\op{Trace}(\E_z=\E_{F^{n}(z)}\to  \dots\to \E_z)} {\det (1-(F^{n})'_z)}
     \right)$$
    is the Taylor expansion of a  rational function in $t$. It seems that
    in many cases for {\it non-compact} varieties $X$ a weaker form of rationality holds as well, when no
    equivariant compactification can be found. Namely, the above
    series (called  the Ruelle zeta-function in general, not necessarily algebraic case) admits a meromorphic continuation to $\C$; also  
    the zeta-function 
    in the
    version  without the denominator is often rational in the non-compact case.

    \subsubsection{Rationality conjecture for motivic local systems}

    In the case hypothetically corresponding to motivic local systems on curves (in the setting of Conjecture 3), one can make a natural a priori
     guess about the denominator in the l.h.s. of the trace formula.
      Namely, for a fixed point $z$ of the map $F^{n}$ corresponding to a fixed point $[\rho]$ in the
       space of representations of $\pi_1(C\times_{\spec \,\F_q}\spec\,\Fb_q)$ in $GL(N,\Qb_l)$, we expect that
        the vector space $T_z X$ together with the automorphism $(F^{n})'_z$ should be
         isomorphic (after the change of scalars) to
         $$H^1(C\times_{\spec \,\F_q}\spec\,\Fb_q,\op{End}(\rho))=\op{Ext}^1(\rho,\rho)$$ endowed
         with the Frobenius operator.

         Eigenvalues of $\Fr^n$ in this case have norm $q^{n/2}$ by the Weil conjecture, hence
         not equal to $1$, and the denominator in the Ruelle  zeta-function does not vanish
          (meaning that the fixed points are non-degenerate).

  In our basic example from Section 0.1 one can propose an explicit formula for the
  denominator term. Define (in notation from Section 0.1)
  for given $t\in\F_q\setminus\{0,1\}$  a matrix
  $T_{tan}\in \op{Mat}(\F_{q}\times \F_{q},\Q)$
 by the formula
 $$T_{tan}:=-\frac{1}{q}\sum_{x\in \F_{q}}
 \left(T_x\right)^2+(q-3-1/q)\cdot id_{\Q^{\F_q}}\,\,\,.$$

 This matrix satisfies the same properties as Hecke operators\footnote{The only difference is that
  eigenvalues of operators $T_x$ are algebraic integers while eigenvalues of $T_{tan}$ are algebraic integers
   divided by $q$.}.
 Namely, all eigenvalues of $T_{tan}$ belong to
 $[-2\sqrt{q},+2\sqrt{q}]$,
 any element of $\spec(T_{tan})$  can be written as $\lambda+\overline{\lambda}$ where $|\lambda|=\sqrt{q}$
is a $q$-Weil number, and for any
$\xi=\lambda+\overline{\lambda}\in \spec(T_x)$ and any integer
$n\ge 1$ the spectrum of the matrix $T_{tan}^{(n)}$ corresponding to
$x\in \F_q\subset\F_{q^n}$ (if we pass
 to the extension $\F_{q^n}\supset \F_q$) contains the element $\xi^{(n)}:=\lambda^n+\overline{\lambda}^n$.

  We expect that the eigenvalue of $T_{tan}$  at the point of the spectrum corresponding to motivic local system $\rho$
  is equal to the trace of
  Frobenius in a two-dimensional submotive of
  the motive $H^1(C,\op{End}(\rho))$, corresponding to the deformations of $\rho$ preserving the unipotency of the monodromy
  around punctures. 
  
  Notice that any motivic local system $\rho$ on $C$ can be  endowed with a non-degenerate skew-symmetric pairing
  with values in the Tate motive. This explains the
  main term of the formula:
  \begin{itemize}
  \item the sum of squares of Hecke operators means that we are using the trace formula
  for Frobenius in the cohomology of $C$ with  coefficients in the tensor square of
  $\rho$,
  \item the factor $1/q$ comes from the Tate twist,
  \item the minus sign comes from the odd (first) cohomology.
  \end{itemize}

  The candidate for the denominator term in the putative Ruelle zeta-function
    is the following operator commuting with the Hecke operators (we write the formula only for the first iteration,
    $n=1$), considered as a function on the spectrum:
    $$D:=\left(q+1-T_{tan}\right)^{-1}\,\,\,.$$
    The reason is that the eigenvalue of $D$ at the eigenvector corresponding to motivic local system
    $\rho$ is equal to
 $$\frac{1}{(1-\lambda)(1-\overline{\lambda})}=\frac{1}{1+q-\xi}
$$
where $\lambda, \overline{\lambda}$ are Weil numbers, eigenvalues
of Frobenius in $H^1(C,\op{End}(\rho))$ satisfying equations
$$\lambda+\overline{\lambda}=\xi
,\,\,\lambda+\overline{\lambda}=q\,\,\,.$$

The  l.h.s. of the putative trace formula for the equivariant vector bundle ${\cal E}_{x_1}\otimes \dots \otimes{\cal E}_{x_k}$
  (here ${\cal E}_x$ is the $F$-equivariant vector bundle corresponding to point $x\in C(\F_q)$, see Conjecture 3), is given 
  (for the $n$-th iteration) by the formula
    $$ \op{Trace} \left(T_{x_1}^{(n)}\dots T_{x_k}^{(n)} D^{(n)}\right)\,\,\,.$$

 It looks that in order to achieve
 the rationality of the putative Ruelle zeta-function one has to add
 by hand certain contributions corresponding to ``missing fixed points''.
   For example, 
     for any $x\in \F_q\setminus\{0,1,t\}$ one has
   $$\op{Trace}(T_x D)=\frac{q}{(q-1)^2}$$
   and the corresponding zeta-function
   $$ \exp\left(-\sum_{n\ge 1} \frac {t^n}{n} \frac{q^n}{(q^n-1)^2}\right)= \prod_{m\ge 1} (1-q^{-m}t)^m\in \Q[[t]]$$
is meromorphic but not rational.
   The  above zeta-function looks
   like the contribution of just one\footnote{Maybe the complete interpretation should be a bit more complicated as
    one can check numerically that $ Trace(D)=\frac{q^2(q-2)}{(q-1)^2(q+1)}$.} fixed point $z_0$ on an
   algebraic dynamical system $z\mapsto F(z)$ on a
   two-dimensional variety, with the spectrum of $(F')_{z_0}$
   equal to $(q,q)$, and the spectrum of the map on the fiber $\E_{z_0}\to
   \E_{z_0}$ equal to $(q,0)$. Here is the precise conjecture coming from computer experiments:

   \begin{conj} For any $x_1,\dots,x_k \in\F_q\setminus\{0,1,t\}, \,\,k\ge 1$
       the series
       $$ \exp\left(
     -\sum_{n\ge 1} \frac{t^n}{n} \left\{
     \op{Trace}\left(T_{x_1}^{(n)}\dots T_{x_k}^{(n)} D^{(n)}\right)+\op{Corr}(n,k)\right\}\right)$$
     where
$$\op{Corr(n,k)}:=-\frac{\left(-1-q^n\right)^k }{\left(1-q^{-n}\right)\left(1-q^{2n}\right)  }\,\,,$$
 is a rational function.
      \end{conj}

The rational function in the above conjecture should be an
L-function of a motive over $\F_q$,  all its zeroes and poles
should be $q$-Weil numbers. 

 Finally, if one considers Ruelle zeta-functions {\it without} the
 denominator term, then rationality is elementary, as will become clear
 in the next section.
 
 \section{Second proposal: formalism of motivic function spaces and higher-dimensional Langlands
correspondence}

 The origin of this section is property 6 (the multiplication table) of Hecke operators in our example from Section 0.1.
\subsection{Motivic functions and the tensor category $\CC_\k$}

Let $S$ be a noetherian scheme.
\begin{dfn} The commutative ring $\FF(S)$ of poor man's motivic
functions\footnote{The name was suggested by V.~Drinfeld.} on $S$ is
the quotient of the free abelian group generated by equivalence
classes of schemes of finite type over $S$, modulo relations
$$[X\to S]=[Z\to S]+[(X\setminus Z)\to S]$$
where $Z$ is a closed subscheme of $X$ over $S$. The
multiplication on $\FF(S)$ is given by the fibered product over
$S$.
\end{dfn}
In the case when $S$ is the spectrum of a field $\k$, we obtain
the standard definition\footnote{Usually one extends the Grothendieck ring of
varieties by inverting the class $[A^1_\k]$ of the affine line,
which is the geometric counterpart of the inversion of
the Lefschetz motive $L=H_2(\P_k^1)$ in the construction
 of Grothendieck pure motives. Here also we can do the same thing.} of the Grothendieck ring of varieties over
 $\k$.
Any motivic function on $S$ gives a function on the set of points of $S$ with values in the
 Grothendieck
rings corresponding to the residue fields.

For a given field $\k$ let us consider the following additive category $\CC_\k$.
Its objects  are schemes of finite type over $\k$, the abelian groups of homomorphisms
are defined by
$$\op{Hom}_{\CC_\k}(X,Y):=\FF(X\times Y)\,\,.$$
The composition of two morphisms represented by  schemes is given by the
fibered product as below:
$$ [B\to Y\times Z]\circ[A\to X\times Y]:=[A\times_Y B\to X\times Z]$$
and extended by additivity to all motivic functions. The identity morphism $id_X$
is given by the diagonal embedding $X\hookrightarrow X\times X$.

One can start from the beginning from constructible sets over $\k$ instead of schemes.
 The category of constructible sets over $\k$ is a full subcategory of $\CC_\k$,
the morphism in  $\CC_\k$ corresponding to a constructible map $f:X\to Y$ is given by
$[X \stackrel{(id_X,f)}{\tto} X\times Y]$, the graph of $f$.

Finite sums (and products) in $\CC_\k$ are given by the disjoint union.

 We endow category $\CC_\k$ with the following tensor structure on objects:
$$X\otimes Y:=X\times Y$$
and by a similar formula on morphisms. The unit object $\1_{\CC_\k}$ is the point $Spec(\k)$.
 The category
$\CC_\k$ is rigid, i.e. for every object $X$ there exists another object
 $X^\vee$ together with morphisms $\delta_X:X\otimes X^\vee\to \1,\,\,\epsilon_X:\1\to X^\vee\otimes X$
 such that both compositions:
$$ X \stackrel{id_X\otimes \epsilon_X} {\tto} X\otimes X^\vee\otimes X \stackrel{\delta_X\otimes id_X} {\tto}
 X,\,\,
X^\vee\stackrel{ \epsilon_X\otimes id_{X^\vee}}{\tto} X^\vee\otimes X\otimes X^\vee
 \stackrel{id_{X^\vee}\otimes \delta_X} {\tto} X^\vee $$
are identity morphisms. The dual object $X^\vee$ in $\CC_\k$ coincides with $X$,
the duality morphisms $\delta_X,\,\,\epsilon_X$ are given by the diagonal embedding
$X\hookrightarrow X^2$.

As in any tensor category, the ring $\op{End}_{\CC_\k}(\1_{\CC_\k})$ is commutative, and the whole
 category is linear over this ring, which is nothing but  the Grothendieck ring of varieties over $\k$.

\subsubsection{Fiber functors for finite fields}
If  $\k=\F_q$ is a finite field then there is an infinite chain $(\phi_n)_{n\ge 1}$
of tensor functors from $\CC_\k$ to the category of finite-dimensional vector spaces over $\Q$.
 It is defined on objects by the formula
$$\phi_n(X):=\Q^{X(\F_{q^n})}\,\,\,.$$
The operator corresponding by $\phi_n$ to a morphism $[A\to X\times Y]$ has the following matrix coefficient with indices
$(x,y)\in X(\F_{q^n}) \times Y(\F_{q^n})$:
$$\#\{a\in A(\F_{q^n}) \,|\,\,a\mapsto (x,y) \}\in \Z_{\ge 0}\subset\Q$$

The functor $\phi_n$ is not canonically defined for $n\ge 2$, the ambiguity is the cyclic group
$\Z/n\Z=\Gal(\F_{q^n}/\F_q)\subset \op{Aut}(\phi_n)$.

\subsubsection{Extensions and variants}

The abelian group $\FF(S)$ of poor man's motivic functions can (and probably should) be replaced
 by the $K_0$ group of the triangulated category ${Mot}_{S, \Q}$ of ``constructible motivic sheaves''
 (with coefficients
 in $\Q$) on $S$.
 Although the latter category is not yet rigorously defined, one can  envision
 a reasonable candidate for the elementary description of $K_0({Mot}_{S,\Q})$. This group should be generated
 by equivalence classes of families of Grothendieck motives (with coefficients in $\Q$) over closed
subschemes of $S$, modulo a suitable equivalence relation. Moreover, group $K_0(Mot_{S,\Q})$ should be filtered
 by the dimension of support, the associated graded group should be canonically isomorphic
to the direct sum over all points $x\in S$ of $K_0$ groups of
categories of pure motives (with coefficients in $\Q$) over the
residue fields\footnote{I do not know how to fill all the details in
the above sketch.}.

Similarly, one can extend the coefficients of the motives from $\Q$ to any field of zero characteristic.
 This change will affect the group $K_0$ and give a different algebra of motivic functions.

Finally, one can add  formally images of projectors to the category $\CC_\k$.

\begin{que} Are there interesting non-trivial projectors in $\CC_\k$?
\end{que}

I do not know at the moment any example of an object in the Karoubi closure of $\CC_\k$
 which is not isomorphic to a scheme. Still, there are interesting non-trivial isomorphisms
 between objects of $\CC_\k$, for example the following version of the Radon transform.

\subsubsection{Example: motivic Radon transform}

Let  $X=\P(V)$ and $Y=\P(V^\vee)$ be two dual projective spaces over $\k$. We assume that $n:=\dim V$ is at
least $3$.

The incidence relation gives a subvariety $Z\subset X\times Y$, which can be interpreted
 as a morphism in $\CC_\k$ in two ways:
$$f_1:=[Z\hookrightarrow X\times Y]\in \op{Hom}_{\CC_\k}(X,Y)\,$$
$$f_2:=[Z\hookrightarrow Y\times X]\in \op{Hom}_{\CC_\k}(Y,X)$$

The composition $f_2\circ f_1$ is equal to
$$[\A^{n-2}] \cdot id_X+ [\P^{n-3}]\cdot [X\to pt \to X]\,\,\,.$$
The reason is that the scheme of hyperplanes passing through points $x_1,x_2\in X$
 is either $\P^{n-3}$ if $x_1\ne x_2$, or  $\P^{n-2}$ if $x_1=x_2$. On the level on constructible sets
 one has
$\P^{n-2}=\P^{n-3}\sqcup \A^{n-2}$.

The first term is the identity morphism multiplied by the $(n-2)$-nd power of the Tate
motive, while the second term is in a sense a rank $1$ operator. It can be killed after passing to the
 quotient of $X$ by $pt$ which is in fact a direct summand in $\CC_\k$:
$$X\simeq pt\oplus (X\setminus pt)\,\,\,.$$
Here we have to choose a point $pt\in X$.
 Similar arguments work for $Y$, and as the result we obtain an isomorphism (inverting the Tate motive)
 $$X\setminus pt\simeq
Y\setminus pt$$
in the category $\CC_\k$ which is not a geometric isomorphism of constructible sets.
\subsection{Commutative algebras in $\CC_\k$}

By definition, a  unital commutative associative algebra $A$ in the tensor category
$\CC_\k$ is given by a scheme of finite type $X/\k$, and  
 two elements
$$\1_A\in \FF(X),\,\,m_A\in \FF(X^3)$$
 (the unit and the product in $A$) satisfying the usual constraints of unitality, commutativity and
 associativity.

The formula for the structure constants $c_{xyz}=(T_x)_{yz}$ of the algebra of Hecke operators
 in our basic example (see 0.1) is given explicitly by counting points
 on constructible sets depending constructibly on a point $(x,y,z)\in X^3$ where $X=\A^1$,
for any $t\in \k\setminus\{0,1\}$ (one should replace factors $q$  by bundles with  fiber $\A^1$).
Hence we obtain  a motivic function on $X^3$ which gives the structure
 of a commutative algebra on $X$ for any $t\in \k\setminus\{0,1\}$,
  for arbitrary field $\k$. 
A straightforward check (see Proposition 1 in Section 2.4 below for a closely related statement)
shows that this algebra is  associative. 

\subsubsection{Elementary examples of algebras}

The first example of a commutative algebra is given by an arbitrary scheme $X$ (or a constructible set)
of finite type over $\k$. The multiplication tensor is given by the diagonal embedding
 $X\hookrightarrow X^3$, the unit is given by the identity map $X\to X$.
 If $\k=\F_q$ is finite then for any $n\ge 1$ the algebra $\phi_n(X)$ is the algebra
 of $\Q$-valued functions
 on the finite set $X(\F_{q^n})$, with the pointwise multiplication.

The next example corresponds to the case when $X$ is
  an abelian group scheme (e.g. $\GG_a$, $\GG_m$, or an abelian variety).
  We define the multiplication tensor $m_A\in \FF(X^3)$
 as the graph of the multiplication morphism $X\times X\to X$. Again, if $\k$ is finite
then the algebra $\phi_n(X)$ is the group algebra with coefficients in $\Q$ of the finite
abelian group $X(\F_{q^n})$. Its points in $\Qb$ are additive (resp. multiplicative) characters
 of $\k$ if $X=\GG_a$ (resp. $X=\GG_m$).

Also, one can see that the algebra in $\CC_\k$ corresponding to the group scheme $\GG_a$ is isomorphic to the direct sum of
  $\1_{\CC_\k}$ (corresponding to the trivial additive character of $\k$) and another algebra $A'$ which can be thought as
 parameterizing {\it non-trivial} additive characters of the field, with the underlying scheme  $\A^1\setminus \{0\}$.
 
 Finally, one can make ``quotients'' of abelian group schemes by finite groups of automorphisms.
 For example, for $\GG_a$ endowed with the action of the antipodal involution $x\to -x$, the formula for the product for
  the corresponding algebra is the sum of the following ``main term''
  $$[Z\mono (\A^1)^3],\,\,\,\,\,Z=\{(x,y,z) |\,\,x^2+y^2+z^2-2(xy+yz+zx)=0\}$$
  (the latter equation means that $\sqrt{x}+\sqrt{y}=\sqrt{z}$), and of certain correction terms.
  Similarly, for the antipodal involution $(x,w)\to (x,-w)$ on the elliptic curve $E\subset \P^1\times \P^1$ given by
   $w^2=x(x-1)(x-t)$ (with $(\infty,\infty)$ serving as zero for the group law), the quotient is $\P^1$ endowed
with the 
multiplication law
    similar to one from the example  0.1. The main term is given by the hypersurface $f_t(x,y,z)=0$ in the notation from Section 0.1.
     The spectrum of the corresponding algebra is rather trivial, in comparison to our example. The difference is that in Section 0.1
      we consider the two-fold cover of $(\A^1)^3$ ramified at the hypersurface $f_t(x,y,z)=0$.
 
 \subsubsection{Categorification}

 One may wonder whether a commutative associative algebra $A$ in $\CC_\k$
 (for general field $\k$, not necessarily finite) is in fact a materialization
of the structure of a symmetric (or only braided)
monoidal category on a triangulated category, i.e.  whether 
 the multiplication morphism is the class in $K_0$ of a bifunctor defining the monoidal structure.
The category under consideration should be either the category of constructible mixed motivic sheaves on the
underlying scheme of $A$,
 or some small modification of it not affecting the group $K_0$
(e.g. both categories could have  semi-orthogonal decompositions with
 the same factors).

\subsection{Algebras parameterizing motivic local systems}

As we noticed already, the Example 0.1 can be interpreted as a commutative associative algebra in $\CC_\k$ 
 parameterizing in a certain sense (via the chain of functors $(\phi_n)_{n\ge 1}$) motivic local systems on a curve over $\k=\F_q$. 
 Here we will formulate a general conjecture, which goes beyond the case of curves.

\subsubsection{Preparations on ramification and motivic local systems}

Let $Y$ be a smooth geometrically connected projective variety over
a finitely generated field $\k$. Let us denote by $K$ the field of
rational functions on $X$ and by $K'$ the field of rational
functions on
 $Y':=Y\times_{\spec\,\k}\spec\,\kb$.
We have an exact sequence
$$1\to \Gal(\overline{K}/K')\to \Gal(\overline{K}/K)\to \Gal(\kb/\k)\to 1$$
For a continuous homomorphism
$$\rho: \Gal(\overline{K}/K')\to GL(N,\Qb_l)$$
where $l\ne \op{char}(\k)$,
which factorizes through the quotient $\pi_1^{geom}(U)$
for some open subscheme
 $U\subset Y'$
one can envision some notion of  ramification divisor (similar to
the notion of the conductor in one-dimensional case)
 which should be  an effective divisor on $Y'$.

One expects that for a pure motive of rank $N$ over $K$
 with coefficients in $\Qb$, the ramification divisor
 of the corresponding $l$-adic local system does not depend on prime $l\ne \op{char}(\k)$, at least for large $l$.

Denote by $\op{IrrRep}_{Y',N,l}$ the set of conjugacy classes of
irreducible representations $\rho: \Gal(\overline{K}/K)\to
GL(N,\Qb_l)$ factorizing through $\pi_1^{geom}(U)$
for some open subscheme
 $U\subset Y'$
 as above.
The Galois group $\Gal(\kb/\k)$ acts on this set.

 Denote by $\op{IrrRep}^{mot,geom}_{Y,N}$ the set of equivalence classes of pure motives in the sense of Grothendieck (defined using
 the numerical equivalence) of rank $N$ over $K$, with
coefficients $\Qb$, which
 are absolutely simple (i.e. remain simple after the pullback to $K'$), modulo the action
 of the Picard group of rank $1$ motives over $\k$ with coefficients in $\Qb$. This set is endowed with
a natural action of $\Gal(\Qb/\Q)$.
  The superscript $geom$ indicates that we are interested only in representations of the geometric
 fundamental group.

One expects that the natural map from $\op{IrrRep}^{mot,geom}_{Y,N}$ to
 the set of fixed points $(\op{IrrRep}_{Y',N,l})^{\Gal(\kb/\k)}$ is a bijection.
In particular, it implies that one can define the ramification divisor for an element of
$\op{IrrRep}^{mot,geom}_{Y,N}$. Presumably, one can give a purely geometric definition of it, without referring
 to $l$-adic representations.

\subsubsection{Conjecture on algebras parameterizing motivic local systems}

\begin{conj} For a smooth projective geometrically connected variety
$Y$  over a finite field $\k=\F_q$, an effective divisor $D$ on $Y$, 
 and a positive integer $N$,  there exists a commutative associative unital algebra $A=A_{Y,D,N}$ in the category $\CC_{\k}$ 
satisfying the following property:

 For any $n\ge 1$ the algebra
$\phi_n(A)$ over $\Q$ is semisimple (i.e. it is a finite direct
sum of number fields) and for any prime $l,\,\,(l,q)=1$ there
exists a bijection between $\op{Hom}_{\Q-alg}(\phi_n(A),\Qb)$ and the
set
 of elements of $\op{IrrRep}^{mot,geom}_{Y\times_{\spec\,\F_q} \spec\,\F_{q^n},N}$ for which the
 ramification divisor is $D$.
 Moreover, the above bijection is
 equivariant with the respect to the natural $Gal(\Qb/\Q)\times \Z/n\Z$-action.
\end{conj}

One can also try to formulate a generalization of the above conjecture, allowing not
 an individual variety $Y$ but a family, i.e. a smooth projective morphism ${\cal{Y}}\to B$ to a
scheme of finite type over $k$, with geometrically connected fibers, together with
a flat family of  ramification divisors.
The corresponding algebra should parameterize choices of a point $b\in B(\F_{q^n})$
 and a irreducible motivic system of given rank and a given ramification on the fiber ${\cal{Y}}_b$.
 This algebra should  map to the  algebra of functions with the pointwise product (see 2.2)
 associated with the base $B$. 
 
In the above conjecture we did not describe how to associate a {\it tower} of finite sets to the algebra
 $A$, as a priori we have just a {\it sequence} of finite sets $X_n:= \op{Hom}_{\Q-alg}(\phi_n(A),\Qb)$
  without no obvious maps between them. This leads to the following
  \begin{que} Which property of an associative commutative algebra $A$ in $\CC_{\F_q}$ gives naturally
   a chain of embeddings
      $$\op{Hom}_{\Q-alg}(\phi_{n_1}(A),\Qb)\mono\op{Hom}_{\Q-alg}(\phi_{n_1 n_2}(A),\Qb)$$
      for all integers $n_1,n_2\ge 1$ ?
  \end{que}
It looks that this holds  automatically, by a kind of trace morphism.

\subsubsection{Arguments in favor, and extensions}
 First of all, there is a good reason to believe that Conjecture 5 holds for curves. Also, it would be reasonable
 to consider local systems with an arbitrary structure group $G$   instead of $GL(N)$.  The algebra parameterizing motivic local system
 on curve $Y=C$ with structure group $G$ should be (roughly) equal to some finite open part of the moduli stack $Bun_{G^L}$ 
  of $G^L$-bundles on $C$, where $G^L$ is the Langlands dual group.
  The multiplication should be given by the class of a motivic constructible sheaf on 
  $$(Bun_{G^L})^3=Bun_{G^L}\times Bun_{G^L\times G^L}$$
  which should be a geometric counterpart  to the lifting of automorphic forms corresponding to the diagonal embedding
  $$G^L\to G^L\times G^L\,\,\,.$$
  Presumably, the multiplication law from Example 0.1 corresponds to the lifting.

If we believe in the Conjecture 5 in the case of curves, then it is
very natural to believe in it in general. The reason is that for a
higher-dimensional variety $Y$ (not necessarily compact) there
exists a curve $C\subset Y$ such that
 $\pi_1^{geom}(Y)$ is a {\it quotient} of $\pi_1^{geom}(C)$. Such a curve
 can be e.g. a complete intersection of ample divisors, the surjectivity
 is a particular case of the Lefschetz theorem on hyperplane sections.
Therefore, the set of equivalence classes of absolutely irreducible
motivic local systems on $Y\times_{\spec\,\F_q}\spec\,\F_{q^n}$ should be
 a {\it subset} of the corresponding set for $C$ for any $n\ge 1$,
 and invariant under $\Gal(\Qb/\Q)$-action as well.
 It looks very plausible that such a collection of subsets should arise
 from a quotient algebra in $\CC_{\F_q}$.

From the previous discussion it looks that the motivic local systems
 in  higher-dimensional
case are ``less interesting'', the $1$-dimensional case is the richest one.
 Nevertheless, there is definitely a  non-trivial higher-dimensional information
 about local systems which can not be reduced to $1$-dimensional data. Namely, for any motivic local system $\rho^{arith}$ on $Y$ and an
 integer $i\ge 0$ 
 the cohomology space
 $$H^i(Y',\rho)$$
 where $\rho$ is the pullback of $\rho^{arith}$ to $Y'$, is a motive over the finite field $\k=\F_q$. We can calculate
 the trace of $N$-th power of Frobenius on it for a given $N\ge 1$, and get a $\Qb$-valued function\footnote{Here there is a small ambiguity
  which should be resolved somehow,
  as one can multiply $\rho^{arith}$ by a one-dimensional motive over $\k$ with coefficients in $\Qb$.} on the set 
 $$X_n:=\op{Hom}_{\Q-alg}(\phi_n(A),\Qb)$$
 for each $n\ge 1$.
This leads to a natural addition to Conjecture 5. Namely, we expect that systems of $\Qb$-valued functions on $X_n$
 associated with higher cohomology spaces arise from elements in
 $\op{Hom}_{\CC_{\F_q}}(\1, A)$ (i.e. from motivic functions on the constructible set underlying algebra $A$).
 
 More generally, one can expect that the motivic constructible sheaves with some kind of boundedness will be parametrized
  by commutative algebras.

Formulas from the example from Section 0.1 make sense and give an algebra in $\CC_\k$
for arbitrary field $\k$. This leads to
\begin{que} Can one construct algebras $A_{Y,D,N}$ for arbitrary ground field $\k$,
 not necessarily finite?
 In what sense will  these algebras ``parameterize'' motivic local system?
\end{que}
  In general, it seems that the natural source of commutative algebras in $\CC_\k$ is not the representation theory, but (quantum)
  algebraic integrable systems. 

\subsection{Towards integrable systems over local fields}

Here we will describe briefly an analog of commutative algebras of
integral operators as above for arbitrary local fields, i.e.  $\R,\,\C$, or finite extensions of
$\Q_p$ or $\F_p((x))$. 
 Let us return to our basic example.
  The check of the associativity of the multiplication law given
  by formula from Section 0.1 in the case of finite fields is reduced to an 
  identification of certain varieties. The most non-trivial part
  is the following
  \begin{prp} For generic parameters $t,x_1,x_2,x_3,x_4$ the two
  elliptic curves
  $$ E:\,\,f_t(x_1,x_2,y)=w_{12}^2,\,\,f_t(y,x_3,x_4)=w_{34}^2 $$
  $$\tilde{E}:\,\,f_t(x_1,x_3,\tilde{y})=\tilde{w}_{13}^2,\,\,f_t(\tilde{y},x_2,x_4)=\tilde{w}_{24}^2 $$
given by equations in variables $(y,w_{12},w_{34})$ and
$(\tilde{y},\tilde{w}_{13},\tilde{w}_{24})$ respectively, are
canonically isomorphic over the ground field. Moreover, one can
choose such an isomorphism which identifies the abelian differentials
 $$ \frac{dy}{w_{12}w_{34}}\mbox{ and } \frac{d
 \tilde{y}}{\tilde{w}_{13},\tilde{w}_{24}}\,\,.$$
    \end{prp}
 In fact, it is enough to check the proposition over an
 algebraically closed field and observe that the curves $E,\tilde{E}$
 have points over the ground field\footnote{Curve $E$ has $16$ rational points with coordinate
$y\in\{0,1,t,\infty\}$, same for $\tilde{E}$. }.

Let now $\k$ be a local field. For a given $t\in
\k\setminus\{0,1\}$ we define a (non-negative) half-density
$c_t$ on $\k^3$ by the formula $$
c_t:=\pi_*\left(\frac{|dx_1|^{1/2}|dx_2|^{1/2}|dx_3|^{1/2}}{|w|}\right)$$
where
$$\pi:Z(\k)\to \A^3(\k),\,\,\,\pi(x_1,x_2,x_3,w)=(x_1,x_2,x_3)$$
is the projection of the hypersurface $$Z\subset\A^4_\k:\,\,
f_t(x_1,x_2,x_3)=w^2\,\,\,.$$ We will interpret $c_t$ as a half-density
on $ (\P^1(\k))^3$ as well.

 One can deduce from the above Proposition the
following
\begin{thm}
The operators $T_x,\,\,x\in \k\setminus\{0,1,t\}$ on the Hilbert space
of half-densities on $\P^1(\k)$, given by
$$T_x(\phi)(y)=\int_{z\in \P^1(\k)} c_t(x,y,z)\,\,\phi(z)$$
 are commuting compact self-adjoint operators.
\end{thm}

Moreover, in the non-archimedean case one can show that the joint
spectrum of commuting operators as above is discrete and consists
of densities locally constant  on $\P^1(\k)\setminus\{0,1,t,\infty\}$. In particular, all eigenvalues of
operators $T_x$ are algebraic complex numbers. Passing to the limit over finite extensions of $\k$ we obtain a countable set upon which acts
  $$ \Gal(\Qb/\Q)\times \Gal(\kb/\k)\,\,\,.$$

Also notice that in the case of local fields the formula is much
simpler then the motivic one, there is no correction terms.
 On the other hand, one has a new ingredient, the local density of
 an integral operator. In general, one can imagine a new
 formalism\footnote{A somewhat similar formalism was proposed by Braverman and Kazhdan (see \cite{BK}, who had in mind
 orbital integrals in the usual local Langlands correspondence.}
 where the structure of an algebra is given by data $(X,Z,\pi,\nu)$
  where $X$ is a (birational type of)  variety over a given field
  $\k$, $Z$ is another variety, $\pi:Z\to X^3$ is a map (defined
  only at the generic point of $Z$), and $\nu$ is a rational section
  of line bundle $K_Z^{\otimes 2}\otimes
  \pi^*(K_{X^3}^{\otimes{-1}})$. If $\k$ is a local field then the
  pushdown by $\pi$ of $|\nu|^{1/2}$ is a half-density on $X^3$.
  The condition of the associativity would follow from a property of  certain data formulated purely in terms of
 birational algebraic
  geometry.

  Presumably, the spectrum for the case of the finite field is
  just a ``low frequency'' part of much larger spectrum for $p$-adic fields,
  corresponding to some mysterious objects\footnote{It looks that all this goes beyond motives, and on the automorphic side
  is related to some kind of  Langlands correspondence for two (or more)-dimensional mixed local-global fields.}.
  
  The commuting integral operators in the archimedean case $\k=\R,\C$ are similar to ones found recently in the usual quantum algebraic 
  integrable systems,
   see \cite{Lebedev}.
\section{Third proposal: lattice models}

\subsection{Traces depending on two indices}

Let $X$ be a constructible set over $\F_q$ and $M$ be an
endomorphism of $X$ in the category $\CC_{\F_q}$ (like e.g. a Hecke operator).  What kind of
object can be called the ``spectrum'' of $M$?

Applying the functors $\phi_n$ for $n\ge 1$  we obtain an infinite
sequence of finite matrices, of exponentially growing size.
We would like to understand  the behavior of spectra of operators
$\phi_n(M)$ as $n\to+\infty$.
 A similar question arises in some models in quantum physics where one is interested in
 the spectrum of a system with finitely many states, with the dimension of the
Hilbert space depending exponentially on the ``number of particles''.

Spectrum of an operator acting on a finite-dimensional space can
be reconstructed from traces of all positive powers. This leads us
to the consideration of the following collection of numbers
$$Z_M(n,m):=\op{Trace}((\phi_n(M))^m)$$
where $n\ge 1$ and $m\ge 0$ are integers. It will be important
later to restrict attention only to strictly positive  values of
$m$, which mean that we are interested only in non-zero
eigenvalues of matrices $\phi_n(M)$, and  want to ignore the
multiplicity of the zero eigenvalue.

{\bf Observation 1.} For a given $n\ge 1$ there exists a finite collection of non-zero
complex numbers $(\lambda_i)$ such that for any $m\ge 1$ one has
$$Z_M(n,m)=\sum_i \lambda_i^m\,\,\,.$$

{\bf Observation 2.} For a given $m\ge 1$ there exists a finite collection of non-zero
complex numbers $(\mu_j)$ and signs $(\epsilon_j\in\{-1,+1\})$, such that for any $n\ge 1$
one has
$$Z_M(n,m)=\sum_j \epsilon_j \mu_j^n\,\,\,.$$

The symmetry between parameters $n$ and $m$ (modulo a minor difference with signs) is quite striking.

The first observation is completely trivial. For a given $n$ the numbers $(\lambda_i)$ are
all non-zero eigenvalues of the matrix $\phi_n(M)$.

 Let us explain the second observation. By functoriality we have
$$Z_M(n,m)=\op{Trace}(\phi_n(M^m))\,\,\,.$$
Let us assume first that $M$ is given by a constructible set $Y$ which maps to $X\times X$:
$$Y\to X\times X,\,\,\,y\mapsto (\pi_1(y),\pi_2(y))\,\,\,.$$
Then $M^m$ is given by the consecutive fibered product
$$Y^{(m)}=Y\times_X Y\times_X\dots\times_X Y\subset Y\times\dots \times Y$$
of $m$ copies of $Y$:
$$Y^{(m)} (\Fb_q)=\{(y_1,\dots,y_m)\in
(Y(\Fb_q))^m| \,\pi_2(y_1)=\pi_1(y_2),\dots,\pi_2(y_{m-1})=\pi_1(y_m)\}$$
The projection to $X\times X$ is given by $(y_1,\dots,y_m)\mapsto (\pi_1(y_1),\pi_2(y_m))$.
To take the trace we should intersect $Y^{(m)}$ with the diagonal. The conclusion is that
$Z_M(n,m)$ is equal to the number of $\F_{q^n}$-points of the constructible set
$$\widetilde{Y}^{(m)}:=Y^{(m)}\times_{X\times X} X\,\,\,,$$
$$\widetilde{Y}^{(m)    }(\Fb_q)=\{(y_1,\dots,y_m) \in Y^{(m)}(\Fb_q)|\, \pi_1(y_1)=\pi_2(y_m)\}\,\,\,.$$
The second observation is now an immediate corollary of the Weil
conjecture on numbers of points of  varieties over finite
fields\footnote{Here we mean only the fact that the zeta-function
of a variety over is rational in $q^s$, and  not the more
deep statement about the norms of Weil numbers.}. The general
case when $M$ is given by a formal {\it integral} linear
combination $\sum_\alpha n_\alpha [Y_\alpha\to X\times X]$ can be treated in a similar way.

\subsection{Two-dimensional translation invariant lattice models}

There is another source of numbers depending on two indices with a
similar behavior with respect to each of indices when another one
is fixed. It comes from the so-called lattice models in statistical
physics.
 A typical example is  the  Ising model.
There is a convenient way to encode Boltzmann weights of a general
lattice model on $\Z^2$ in terms of linear algebra.

\begin{dfn} Boltzmann weights of a  $2$-dimensional translation
 invariant lattice model
are given by a pair $V_1,V_2$ of finite-dimensional vector spaces over $\C$ and a linear operator
 $$R: V_1\otimes V_2\to V_1\otimes V_2\,\,\,.$$
\end{dfn}

Such  data give a function (called the partition function) on a
certain set of graphs. Namely, let $\G$ be a finite oriented graph
whose edges are colored by $\{1,2\}$
 in such a way that for every vertex $v$ there are exactly two edges colored by $1$ and $2$ with head  $v$,
 and also there are exactly two edges colored by $1$ and $2$ with tail  $v$.
Consider the tensor product of copies of $R$ labelled by the set
$Vert(\G)$ of vertices of $\G$. It is an element $v_{R,\G}$ of the
vector space
$$(V_1^\vee\otimes V_2^\vee\otimes V_1\otimes V_2)^{\otimes Vert(\G)}\,\,\,.$$
The structure of an oriented colored graph gives an identification of the above space with
$$(V_1\otimes V_1^\vee)^{\otimes Edge_1(\G)}\otimes (V_2\otimes V_2^\vee)^{\otimes Edge_2(\G)}$$
where $Edge_1(\G), Edge_2(\G)$ are the sets of edges of $\G$ colored
by $1$ and by $2$.
 The tensor product of copies of the standard pairing gives a linear functional $u_\G$
on the above space.
We define the {\it partition function} of the lattice model on $\G$ as
$$Z_R(\G)=u_\G(v_{R,\G})\in \C\,\,\,.$$

An oriented colored  graph $\G$  as above is the same as a finite
set with two permutations $\tau_1,\tau_2$. The set here is
$Vert(\G)$, and permutations $\tau_1,\tau_2$ correspond to edges
colored by $1$ and $2$ respectively.

In the setting of {\it translation invariant} $2$-dimensional
lattice models we are interested
 in the values of the partition function only on  graphs corresponding to pairs of
commuting permutations. Such a graph (if it is non-empty and
connected) corresponds
 to a subgroup $\Lambda\subset \Z^2$ of finite index. We will denote the partition
function\footnote{In physical literature it is called
 the partition function with periodic boundary conditions.} of the graph
 corresponding to $\Lambda$ by $Z_R^{lat}(\Lambda)$.

Finally, Boltzmann data make sense in an arbitrary rigid tensor category $\CC$. The partition
function of a graph
 takes values in the commutative ring $\op{End}_\CC(\1)$.
 In particular, one can speak about {\it super} Boltzmann data for the category $Super_\C$
of finite-dimensional complex
 super vector spaces.

 \subsubsection{Transfer matrices}

Let us consider a special class of lattices $\Lambda\subset \Z^2$ depending on two parameters.
Namely, we set
 $$\Lambda_{n,m}:=\Z\cdot(n,0)\oplus \Z\cdot (0,m)\subset\Z^2\,\,\,.$$
\begin{prp}
For any Boltzmann data $(V_1,V_2,R)$ and a given $n\ge 1$ there exists a finite collection of non-zero
complex numbers $(\lambda_i)$ such that for any $m\ge 1$ one has
$$Z_R^{lat}(\Lambda_{n,m})=\sum_i \lambda_i^m\,\,\,.$$
\end{prp}

The proof is the following. Let us  introduce a linear operator
(called the {\it transfer matrix}) by formula:
 $$T_{(2),n}:=\op{Trace}_{V_1^{\otimes n}} ((\sigma_n\otimes  id_{V_2^{\otimes n}})\circ R^{\otimes n})\in \op{End}(V_2^{\otimes n})$$
 where $\sigma_n\in \op{End}(V_1^{\otimes n})$ is the cyclic permutation. 
Here we interpret $R^{\otimes n}$ as an element of
$$(V_1^\vee)^{\otimes n}\otimes (V_2^\vee)^{\otimes n}\otimes V_1^{\otimes n}\otimes V_2^{\otimes n}=\op{End}(V_1^{\otimes n})\otimes
 \op{End}(V_2^{\otimes n})\,\,\,.$$
It follows from the definition of the partition function that
  $$Z_R^{lat}(\Lambda_{n,m})=\op{Trace}\left(T_{(2),n}\right)^m$$
 for all $m\ge 1$. The collection $(\lambda_i)$ is just the collection of all {\it non-zero}
eigenvalues of $T_{(2),n}$ taken with multiplicities.

Similarly, one can define transfer matrices $T_{(1),m}$ such that
$Z_R^{lat}(\Lambda_{n,m})=\op{Trace}\left(T_{(1),m}\right)^n$ for all $n,m\ge 1$.
We see that the function $(n,m)\mapsto Z_R^{lat}(\Lambda_{n,m})$ has the same two properties as 
 the function $(n,m)\mapsto Z_M(n,m)$ from Section 3.1. For super Boltzmann data one
 obtains sums of exponents with signs.

\subsection{Two-dimensional Weil conjecture}

Let us return to the case of an endomorphism $M\in \op{End}_{\CC_{\F_q}}(X)$. In Section 3.1
 we have defined numbers $Z_M(n,m)$ for $n,m\ge 1$. Results of  3.2 indicate
 that one should interpret pairs $(n,m)$ as parameters for a special class
 of ``rectangular'' lattices in $\Z^2$. A general lattice $\Lambda \subset \Z^2$ depends on
 3 integer parameters
 $$\Lambda=\Lambda_{n,m,k}=\Z\cdot (n,0)\oplus \Z\cdot (k,m),\,\,n,m\ge 1,\,\,0\le k< n\,\,\,.$$
Here we propose an extension of function $Z_M$ to all lattices in $\Z^2$:
$$Z_M(\Lambda_{n,m,k}):=\op{Trace}((\phi_n(M))^m(\phi_n(\Fr_X))^k)$$
where $\Fr_X\in \op{End}_{\CC_{\F_q}}(X)$
 is the graph of the Frobenius endomorphism of the scheme $X$. Notice that $\phi_n(\Fr_X)$
is periodic with period $n$ for any $n\ge 1$.

\begin{prp}
 Function $Z_M$ on lattices in $\Z^2$ defined as above,
 satisfy the following property:
 for any two vectors $\gamma_1,\gamma_2\in\Z^2$
 such that $\gamma_1\wedge\gamma_2\ne 0$ there exists a finite collection of
 non-zero complex numbers $(\lambda_i)$ and signs $(\epsilon_i)$ such that for any $n\ge 1$ one has
$$Z_M(\Z\cdot \gamma_1\oplus \Z\cdot n\gamma_2)=\sum_i \epsilon_i \lambda_i^n\,\,\,.$$
In other words, the series in formal variable $t$
$$\exp\left(-\sum_{n\ge 1} Z_M(\Z\cdot \gamma_1\oplus \Z\cdot n\gamma_2) \cdot t^n/n\right)$$
is rational.
\end{prp}

The proof is omitted here, we'll just indicate that it follows from the consideration of the
 action of the Frobenius element and of cyclic permutations on the (\'etale) cohomology of
 spaces $\widetilde{Y}^{(m)}$
 introduced in Section 3.1.

Also, it is easy to see that the same property holds for the partition function
$Z_R^{lat}(\Lambda_{m,n,k})$ for arbitrary (super) lattice models.\footnote{In general, one can show that for any lattice model given by
operator $R$, and for any  matrix $A\in GL(2,\Z)$ there exists another  lattice model with operator $R'$ such that
 for any lattice $\Lambda\subset \Z^2$ one has $Z_R^{lat}(\Lambda)=Z_{R'}^{lat}(A(\Lambda))$.}
The analogy leads to a two-dimensional analogue of the Weil conjecture
(the name will be explained in the next section):
\begin{conj} For any endomorphism $M\in \op{End}_{\CC_{\F_q}}(X)$ there exists super Boltzmann data
$(V_1,V_2,R)$ such that for any $\Lambda\subset \Z^2$ of finite index one has
$$Z_M(\Lambda)=Z_R^{lat}(\Lambda)\,\,\,.$$
  \end{conj}

Up to now there is no hard evidence for this conjecture, there are just a few
cases where one can construct
 a corresponding lattice model in an ad hoc manner. For example, it is possible (and not totally trivial) to do that for
 the case when $X=\A^1_{\F_q}$ and $M$ is the graph of the map $x\to x^c$ where $c\ge 1$ is an integer.

The above conjecture means that one can see matrices $\phi_n(M)$ as analogs of
transfer matrices\footnote{At least if one is interested in the non-zero part of  spectra.
 In general, the size of the transfer matrix depends on $n$ as an exact exponent, while the size
of $\phi_n(M)$ is a finite alternating sum of exponents.}.
 In the theory of integrable models people are interested in systems where the Boltzmann weights
$R$
 depends non-trivially on a parameter $\rho$
(spaces $V_1,V_2$ do not vary), and the horizontal transfer
 matrices commute with each other
$$[T_{(2),n}(\rho_1),T_{(2),n}(\rho_2)]=0$$
 because of Yang-Baxter equation. Theory of automorphic forms seems to produce
families of commuting endomorphisms in category $\CC_{\F_q}$, which is quite analogous
 to the integrability in lattice models. There are still serious differences.
 First of all,  commuting operators in the automorphic forms case depend on discrete parameters
 whereas in the integrable model case they depend algebraically on continuous parameters.
 Secondly, the spectrum of a Hecke operator in its $n$-th incarnation (like  $T_x^{(n)}$ in Section 0.1)
 has typically $n$-fold degeneracy, which does not happen in
 the case of the usual integrable models with period $n$.

\subsection{Higher-dimensional lattice models and a higher-dimensional
 Weil conjecture} Let $d\ge 0$ be an integer.
\begin{dfn} Boltzmann data of a $d$-dimensional translation invariant lattice model are given by a collection
$V_1,\dots,V_d$ of finite-dimensional vector spaces over $\C$ and a linear operator
 $$R: V_1\otimes \dots \otimes  V_d \to V_1\otimes \dots \otimes V_d\,\,\,.$$
\end{dfn}
Similarly, one can define $d$-dimensional lattice model in an arbitrary rigid tensor category.
The partition function is a function on finite sets endowed with the action of 
the free group with $d$ generators. In particular, for abelian actions, it gives a function
 $\Lambda\mapsto Z_R^{lat}(\Lambda)\in\C$ on the set of subgroups of finite index in $\Z^d$.
Also, for any lattice $\Lambda_{d-1}\subset \Z^d$
of rank $(d-1)$ and a vector $\gamma\in\Z^d$  such that $\gamma\notin \Q\otimes \Lambda_{d-1}$, the function
 $$n\ge 1\mapsto Z_R^{lat}(\Lambda_{d-1}\oplus \Z\cdot n\gamma)$$
is a finite sum of exponents.
Analogously, for any $d$-dimensional lattice model $R$ and any integer $n\ge 1$ there exists
 its dimensional reduction,  periodic with period
 $n$ in $d$-th coordinate, which is a $(d-1)$-dimensional lattice model ${R}_{(n)}$ satisfying the property
  $$Z_{{R}_{(n)}}(\Lambda_{d-1})=Z_{{R}}(\Lambda_{d-1} \oplus \Z\cdot n\,e_{d})\,,\,\,\,\,\forall 
  \Lambda_{d-1}\subset \Z^{d-1}$$ 
   where $e_{d}=(0,\dots,0,1)\in \Z^{d}=\Z^{d-1}\oplus \Z$ is the last standard basis vector.

\begin{conj} For any $(d-1)$-dimensional lattice model $(X_1,\dots,X_{d-1},M),\,d\ge 1$
in the category $\CC_{\F_q}$, there exists a $d$-dimensional super
lattice model
 $(V_1,\dots,V_{d},R)$ in $Super_\C$ such that for any $n\ge 1$ the numerical
 $(d-1)$-dimensional  model $\phi_n(M)$ gives the same partition function on the set of subgroups
of finite index in $\Z^{d-1}$ as the  dimensional reduction
  $R_{(n)}$.
\end{conj}

In the case $d=1$ this conjecture follows from the usual Weil
conjecture. Namely,   a $0$-dimensional Boltzmann data in $\CC_{\F_k}$
is just an element
$$M\in \op{End}_{\CC_{\F_q}}(\1)=\op{End}_{\CC_{\F_q}}(\otimes_{i\in\emptyset} X_i)$$
of
the Grothendieck group of varieties over $\F_k$ (or of $K_0$ of the category
 of pure motives over $\F_k$ with rational coefficients).
The corresponding numerical lattice models $\phi_n(M)$ are just numbers,
  counting $\F_{q^n}$-points in $M$.
By the usual Weil conjecture these numbers are traces of powers of an operator
 in a super vector space, i.e. values of the partition function for $1$-dimensional
 super lattice model.

Similarly, for $d=2$ one gets the $2$-dimensional Weil conjecture from the previous section.

\subsubsection{Evidence: $p$-adic Banach lattice models}

Let $K$ be a complete non-archimedean field (e.g. a finite
extension of $\Q_p$). We define a  $d$-dimensional
{\it contracting} Banach lattice model
 as follows. The Boltzmann  data consists of
 \begin{itemize}
 \item $2d$ countable generated $K$-Banach spaces  $V_1^{in},\dots,V_d^{in},V_1^{out},\dots,V_d^{out}$,
 \item a bounded   operator $R^{vertices}: V_1^{in}\widehat{\otimes}\dots \widehat{\otimes} V_d^{in}\to
 V_1^{out}\widehat{\otimes}\dots\widehat{\otimes} V_d^{out}\,\,,$
 \item a collection of compact operators
  $R_i^{edges}:V_i^{out}\to V_i^{in},\,\,i=1,\dots,d$.
  \end{itemize}
  Such  data again give a function on oriented graphs with colored edges, in the definition
   one should insert operator $R_i^{edges}$ for each edge colored by index $i,\,\,i=1,\dots,d$.
   In the case of 
  {\it finite-dimensional} spaces $(V_i^{in},V_i^{out})_{i=1,\dots,d}$ we obtain the same partition function as for
  a usual
   finite-dimensional lattice model. Namely, one can set
  $$R:=\left(\otimes_{i=1}^d R_i^{edges}\right)\circ R^{vertices},\,\,\,V_i=V_i^{in},\,\forall i=1,\dots,d$$
  or, alternatively,
   $$\tilde{R}:=R^{vertices}\circ \left(\otimes_{i=1}^d
   R_i^{edges}\right),\,\,\,\tilde{V}_i:=V_i^{out},\,\forall i=1,\dots,d\,\,\,.$$

 In particular, for any contracting Banach model one get a function
$\Lambda\mapsto Z_R^{lat}(\Lambda)\in K$ on the set of sublattices
of $\Z^d$. This function satisfies the property that for any
lattice $\Lambda_{d-1}\subset \Z^d$ of rank $(d-1)$ and a vector
$\gamma\in\Z^d$  such that $\gamma\notin \Q\otimes \Lambda_{d-1}$,
one has
 $$Z_R^{lat}(\Lambda_{n-1}\oplus \Z\cdot n\gamma)=\sum_\alpha \lambda_\alpha^n,\,\,\,\forall n\ge 1$$
  where $(\lambda_\alpha)$ is a
(possibly) countable $Gal(\overline{K}/K)$-invariant collection of
non-zero numbers in $\overline{K}$  (eigenvalues of transfer
operators)
 whose norms tend to zero. Similarly, one can define super Banach
 lattice models.

 Here we announce a result supporting Conjecture 7, the proof is a straightforward extension of the 
Dwork method for the proving of the rationality
  of zeta-function of a variety over a finite field.
  \begin{thm} The Conjecture 7 holds if one allows contracting Banach super lattice models over  a finite extension of $\Q_p$ where
  $p$ is the characteristic of the finite field $\F_q$.
  \end{thm}

\subsection{Tensor category $\AA$ and the Master Conjecture}

Let us consider the following rigid tensor category $\AA$.
 Objects of $\AA$ are finite-dimensional  vector spaces over $\C$.
 The set of morphisms $\op{Hom}_\AA(V_1,V_2)$ is defined as the group $K_0$ of the category
 of finite-dimensional representations of the free (tensor) algebra
$$T(V_1\otimes V_2^\vee):=\bigoplus_{n\ge 0} (V_1\otimes V_2^\vee)^{\otimes n}\,\,\,.$$
A representation of the free algebra by operators in a vector
space $U$ is the same as an action of its generators on $U$, i.e.
a linear map
  $$V_1\otimes V_2^\vee\otimes U \to U\,\,.$$
Using  duality we interpret it as a map
 $$V_1\otimes U \to V_2 \otimes U\,\,.$$
The composition of morphisms is defined by the following formula
on generators:
$$[V_1\otimes U \to V_2 \otimes U]\circ [V_2\otimes U' \to V_3 \otimes U']$$
is equal to
 $$[V_1\otimes(U\otimes U')\to V_3\otimes (U\otimes U')]$$
where the expression in the bracket
 is the obvious composition of linear  maps
$$V_1\otimes U\otimes U'\to V_2\otimes U\otimes U'\to V_3\otimes U\otimes U'\,\,\,.$$

The tensor product in $\AA$ coincides on objects with the tensor
product in $Vect_\C$, the same for the duality.
  The formula for the tensor product on morphisms is an obvious one, we leave details to the reader.

Like in Section 2.1 (Question 3), we can ask the following
\begin{que} Are there interesting non-trivial projectors in $\AA$?\footnote{A similar question about commuting endomorphisms in $\AA$ is 
almost
 equivalent to the study of finite-dimensional solutions of the Yang-Baxter equation.}
\end{que}

We denote by $\AA^{kar}$ the Karoubi closure of $\AA$.

There exists an infinite chain of tensor functors $(\phi_n^\AA)_{n\ge 1}$
from $\AA$ to the category
 of finite-dimensional  vector spaces over $\C$ given  by
$$\phi_n^\AA (V):=V^{\otimes n}$$
 on objects, and by
$$[f:V_1\otimes U\to V_2\otimes U] \stackrel{\phi_n}{\longmapsto}
 \op{Trace}_{U^{\otimes n}}((\sigma_n\otimes id_{V_2^{\otimes n}}f^{\otimes n})\in
\op{Hom}_{Vect_\C} (V_1^{\otimes n},V_2^{\otimes n})$$ on morphisms, where $\sigma_n:U^{\otimes n}\to U^{\otimes n}$ is the 
cyclic permutation.
The cyclic group $\Z/n \Z$ acts by automorphisms of $\phi_n^\AA$.
 Moreover, the generator of the cyclic group acting on $V^{\otimes n}=\phi_n^\AA (V)$
is the image under $\phi^\AA_n$
 of a certain central element $\Fr_V$ in the algebra of endomorphisms
 $\op{End}_\AA(V)$. This ``Frobenius'' element is represented by
the linear map $\sigma:V\otimes U\to V\otimes U$ where $U:=V$ and $\sigma=\sigma_2$
is the permutation. As in the case of $\CC_{\F_q}$, for any $V$ the operator
 $\phi_n^\AA (\Fr_V)$ is periodic with period $n$.

Let us introduce a small modification $\AA'$ of the tensor category
$\AA$. Namely, it will have the same objects (finite-dimensional
vector spaces over $\C$), the morphism groups will be the quotients
$$\op{Hom}_{\AA'}(V_1,V_2):=K_0(T(V_1\otimes V_2^\vee)-\op{mod})/\Z\cdot [triv]$$
where $triv$ is the trivial one-dimensional representation of $T(V_1\otimes V_2^\vee)$ given by zero map
 $$V_1\otimes \1\stackrel{0}{\to} V_2\otimes \1$$
All the previous considerations extend to the case of $\AA'$.

Amazing similarities between  categories $\CC_{\F_q}$ and $\AA'$ suggests the following
\begin{conj} For any prime $p$ there exists a tensor functor
 $$ \Phi_p:  \CC_{\F_p}\to {\AA'}^{kar}$$
 and a sequence of isomorphisms of tensor functors from
 $\CC_{\F_p}$ to $Vect_\C$ for all $n\ge 1$
$$ iso_{n,p}: \phi_n^{\AA'}\circ \Phi_p \simeq i_{Vect_\Q\to Vect_\C}\circ \phi_n$$
where $i_{Vect_\Q\to Vect_\C}$ is the obvious embedding functor from the category of vector spaces
 over $\Q$ to the one over $\C$.
 Moreover, for any $X\in \CC_{\F_p}$ the functor $\Phi_p$ maps the Frobenius element $\Fr_X \in \op{End}_{\CC_{\F_p}}(X)$
 to $\Fr_{\Phi_p(V)}$.
\end{conj}

This conjecture we call the Master Conjecture because it implies
 simultaneously {\it all} higher-dimensional versions of the Weil conjecture at
once, as one has the bijection (essentially by definition)
$$\{\mbox{$(d-1)$-dimensional super lattice models in }\AA'\}\simeq $$
$$\simeq
\{\mbox{$d$-dimensional super lattice models in } Vect_\C\}\,\,\,.$$

\begin{rmk} One can consider a larger category $\AA^{super}$ adding
to objects of $\AA$ super vector spaces  as well.
 The group  $K_0$ in the super case should be defined as the naive $K_0$ modulo the relation
  $$[V_1\otimes U \to V_2 \otimes U]=-[V_1\otimes \Pi (U ) \to V_2 \otimes \Pi(U)]$$
where $\Pi$ is the parity changing functor.
\end{rmk}

 It suffices to verify the Master Conjecture only on the full
symmetric monoidal subcategory of  $\CC_{\F_p}$ consisting of
powers $\left(\A^n_{\F_p}\right)_{n\ge 0}$ of the affine
line. The reason is that
 any scheme of finite type can be embedded
 (by a constructible map) in an affine space $\A^n_{\F_p}$, and the characteristic function of the image of such
an embedding
 as an idempotent in $\op{End}_{\CC_{\F_p}} (\A^n_{\F_p})$.

\subsubsection{Machine modelling finite fields}
 Let us fix a prime $p$.
The object $A:=\A^1_{\F_p}$ of $\CC_{\F_p}$ is a commutative algebra (as well as any scheme of finite type, 
see 2.2.1), with
the product given by the diagonal in its cube.
 The category $\op{Aff}(\CC_{\F_p})$ of ``affine schemes'' in $\CC_{\F_p}$ (i.e. the category opposite
 to the category of commutative associative unital algebras in $\CC_{\F_p}$)
 is closed under finite products. In particular, it makes sense to speak about group-like etc. objects
 in $\op{Aff}(\CC_{\F_p})$. Affine line $A$ is a commutative ring-like object
 in $\op{Aff}(\CC_{\F_p})$, with the operations of addition and multiplication corresponding
 to the graphs of the usual addition and multiplication on $\A^1_{\F_p}$.
 In plain terms, this means that besides the
commutative algebra  structure on $A$
$$m:A\otimes A\to A$$
 we have two coproducts (for the addition and for the multiplication)
$$co-a:A\to A\otimes A,\,\,co-m:A\to A\otimes A$$
which are homomorphisms of algebras, and satisfy the usual bunch
of rules for commutative associative rings,
 including the distributivity law.

If the Master Conjecture 8 is true then it gives an object $V_p:=\Phi_p(A) \in {\AA'}^{kar}$,
 with one product and two coproducts.
One can expect that it is just $\C^p$ as a vector space. For any $n\ge 1$
 the $\AA'$-product on $V_p$ defines
 a commutative algebra structure on $V_p^{\otimes n}$. Its spectrum should be a finite set
 consisting of $p^n$ elements. Two coproducts give operations of addition and multiplication on this set, and we
will obtain a {\it canonical} construction\footnote{Compare with question 2 
 in Section 1.3, and remarks afterwards.} of the finite field $\F_{p^n}$ uniformly for all $n\ge 1$.

Even in the case $p=2$ the construction of such $V_p$ is a formidable task: one should find
3 finite-dimensional super representations of the free algebra in 8 generators, satisfying
 9 identities in various $K_0$ groups.

 \subsection{Corollaries of the Master Conjecture}

 \subsubsection{Good sign: Bombieri-Dwork bound}
 
 One can deduce easily from the Master Conjecture that  for any given $p$ and any system of equations in arbitrary number of variables
  $(x_i)$ where each of equations is of an elementary form like
  $x_{i_1}+x_{i_2}=x_{i_3}$, or $x_{i_1}x_{i_2}=x_{i_3}$ or $x_i=1$,
    the number of solutions of this system over $\F_{p^n}$ is an alternating sum of exponents in $n$, 
with the total
    number of terms bounded by $C^N$ where $C=C_p$ is a constant depending on $p$, and $N$ is the 
number of equations.
    In fact, it is a well-known Bombieri-Dwork bound (and $C$ is an absolute constant\footnote{A straightforward application of \cite{Bombieri} gives
    the upper bound $C\le 17^4$ which is presumably very
    far from the optimal one.}), see \cite{Bombieri}.

\subsubsection{Bad sign: cohomology theories for  motives over finite fields}
Any machine modelling finite field should be defined over a finitely generated commutative ring.
In particular, there should be a machine defined over a number field $K_p$ depending only on the characteristic $p$.
 A little thinking shows that the enumeration of the number of solutions of any given system of equations in 
the elementary form as above, 
  will be expressed as a super trace of an operator
 in a finite-dimensional super vector space defined over $K_p$.
 On the other hand, it looks very plausible that the category of
 motives over any finite field $\F_q$ does not have any fiber
 functor defined over a number field, see \cite{Milne} for a discussion. I think that this
 is a strong sign indicating that
  the Master Conjecture is just wrong!

\section{Categorical afterthoughts}

\subsection{Decategorifications of 2-categories}

Two categories, $\CC_\k$ and $\AA$ introduced in this paper have a
common feature which is also shared (almost) by the category of
Grothendieck motives. The general framework is the following.

Let $\B$ be a $2$-category such that for any two objects $X,Y\in
\B$ the category of $1$-morphisms ${Hom}_\B(X,Y)$ is a small {\it additive}
category, and the composition of $1$-morphisms is a
bi-additive functor. In practice we may ask for categories
${Hom}_\B(X,Y)$ to be triangulated categories (enriched in their
turn by spectra, or by complexes of vector spaces). Moreover, the
composition could be
 only a weak functor (e.g. $A_\infty$-functor), and the
 associativity of the composition could hold only up to (fixed)
 homotopies and higher homotopies. The rough idea is that objects
 of $\B$ are ``spaces'' (non-linear in general), whereas objects
  of the category ${Hom}_\B(X,Y)$ are linear things on the ``product''
  $X\times Y$ interpreted as kernels of some additive functors
  transforming some kind of sheaves from $X$ to $Y$, by taking the pullback from $X$, the tensor product with the kernel on $X\times Y$, and then the direct
   image with compact supports to $Y$.

In such a situation one can define a new ($1$-)category $K^{tr}(\B)$
which is in fact a triangulated category. This category will be
called the {\it decategorification} of $\B$.

The first step is to define a new $1$-category $K(\B)$ enriched by
spectra. It has the same objects as $\B$, the morphism spectrum
 ${Hom}_{K(\B)}(X,Y)$ is defined as the spectrum of $K$-theory of
 the triangulated category $Hom_\B(X,Y)$\footnote{It is well-known that in order to define
  a correct $K$-theory one needs either an appropriate enrichment on $Hom_\B(X,Y)$,
  or a model structure in the sense of Quillen, see e.g. \cite{Schlichting}.}.

   The second step is to make a formal triangulated envelope of
   this category. This step needs nothing, it can be performed for
   arbitrary category enriched by spectra. Objects of the new
   category are finite extensions of formal shifts of the objects
   of $K(\B)$, like e.g. twisted complexes by Bondal and Kapranov.

 At the third step one adds formally direct summands for
 projectors. The resulting category $K^{tr}(\B)$ is the same as the
 full category of compact objects in the category of exact
 functors from $K(\B)^{opp}$ to the triangulated category of
 spectra (enriched by itself).

 Finally, one can define a more elementary pre-additive\footnote{Enriched by abelian groups 
in the plain sense (without higher
 homotopies).} category $K_0(\B)$ by defining $Hom_{K_0(\B)}(X,Y)$
 to be $K_0$ group of triangulated category $Hom_\B(X,Y)$.
  Then we add formally to it finite sums and images of
  projectors.
  The resulting additive Karoubi-closed category will be denoted
  by $K_0^{kar}(\B)$ and called $K_0$-decategorification of $\B$. 
 In what follows we list several examples of decategorifications.
 \subsubsection{Non-commutative stable homotopy theory}
 R.~Meyer and R.~Nest introduced in \cite{NM} a non-commutative analog of the triangulated category of spectra.
  Objects of their category are not necessarily unital $C^*$-algebras, the morphism group from $A$ to $B$
 is defined
 as the bivariant Kasparov theory $KK(A,B)$. One of main observations in \cite{NM} is that this gives a structure
 of a triangulated category on $C^*$-algebras. Obviously this construction has a flavor of the $K_0$-decategorification.
\subsubsection{Elementary algebraic model of bivariant K-theory}
One can define a toy algebraic model of the construction by  Meyer and Nest. For a given base field $\k$
 consider the pre-additive category whose objects are unital associative $\k$-algebras, and the group
 of morphisms from $A$ to $B$ is defined as $K_0$ of the exact category consisting of bimodules 
 ($A^{op}\otimes B$-modules) which are projective and finitely generated as $B$-modules.
 This is obviously a $K_0$-decategorification of a 2-category.
\subsubsection{Non-commutative pure and mixed motives}
  Let us consider the quotient of the category of Grothendieck
  Chow
  motives $Mot_{\k,\Q}$ over given field $\k$ with rational coefficients, 
by an autoequivalence given by the invertible functor
$\Q(1)\otimes\cdot$.  The set of morphisms in this category
 between motives of two smooth projective schemes $X,Y$ is given
 by
 $$ \op{Hom}_{Mot_{\k,\Q}/\Z^{\Q(1)\otimes\cdot}}(X,Y)=\bigoplus_{n\in \Z} \op{Hom}_{Mot_{\k,\Q}}
 (X,\Q(n)\otimes Y)=$$
 $$ = \left(\Q\otimes_\Z \bigoplus_{n\in\Z}\op{Cycles}_{n}(X\times Y)\right)/(\mbox{ rational equivalence })
 = $$
 $$=\Q\otimes_\Z \bigoplus_{n\in\Z}CH^n(X\times Y)=\Q\otimes_\Z K^0(X\times Y)$$
 because the Chern character gives an isomorphism modulo torsion
 between the sum of all Chow groups and $K^0(X)=K_0(D^b(Coh\,X))$,
 the $K_0$ group of the bounded derived category $D(X):=D^b(Coh\,X)$ of coherent
 sheaves on $X$. Finally, the category $D(X\times Y)$ can be interpreted
 as the category of functors $D(Y)\to D(X)$.

 Triangulated categories of type $D(X)$ where $X$ is a smooth projective variety over $\k$ belong to a larger 
class of {\it smooth proper} triangulated $\k$-linear dg-categories (another name is ``saturated categories''), see
 e.g. \cite{KS},\cite{TV}.
    We see that the above quotient category of pure motives is a full subcategory
 of $K_0$-decategorification (with $\Q$ coefficients) of the 
$2$-category of smooth proper $\k$-linear dg-categories. This construction was described recently
 (without mentioning the relation to motives) in \cite{tabuada}.

 Analogously, if one takes the quotient of the Voevodsky triangulated category of mixed motives
 by the endofunctor $\Q(1)[2]\otimes \cdot$,  the resulting triangulated category seems
 to be similar to a full subcategory of the full decategorification of the 
$2$-category of smooth proper $\k$-linear dg-categories.

\subsubsection{Motivic integral operators}
We mentioned already in Section 2.1 that the category $\CC_\k$ should be considered as a $K_0$-decategorification
 of a 2-category of motivic sheaves. A similar $2$-category was considered in \cite{G}
 in the relation to questions in integral geometry and calculus of integral operators with holonomic kernels.
 
\subsubsection{Correspondences for free algebras}
The category $\cal A$ is a $K_0$-decategorification by definition.
\subsection{Trace of an exchange morphism}

Let $G_1,G_2$ be two endofunctors of a triangulated category $\CC$, and
an exchange morphism (a natural transformation)
$$\alpha:G_1\circ G_2\to G_2\circ G_1$$
is given\footnote{We do not assume that $\alpha$ is an
isomorphism.}. Under the appropriate
finiteness condition (e.g. when $\CC$ is smooth and proper) one can define the {\it trace} of $\alpha$, which can be calculated
 in two ways, as the trace of endomorphism of $\op{Tor}(G_1,id_\CC)$ associated with $G_2$ and $\alpha$,
 and as a similar trace with exchanged $G_1$ and $G_2$ (see \cite{kapranov} for a related stuff).
 Passing to powers and natural exchange morphisms constructed from $nm$ copies of $\alpha$:
 $$\alpha_{(n,m)}:G_1^{n}\circ G_2^{m}\to G_2^{ m} \circ G_1^{n}$$ we obtain
 a collection of numbers $Z_{\alpha}(n,m):=\op{Trace}(\alpha_{(n,m)})$ for $n,m\ge 1$.
 It is easy to see that these numbers come from a $2$-dimensional super lattice model.

 Let $\CC=D(X)$ for smooth projective $X$, and functors are given by $F^*$ and by $\E\otimes \cdot$
 where $F:X\to X$ is a map, and $\E$ is a vector bundle endowed with a morphism $g:F^*\E\to\E$
 (as in Section 1.3). In this case $Z_\alpha(n,m)$ is the trace (without the denominator)
associated with the map $F^{ n}$ and
 the bundle $\E^{\otimes m}$. For example,
 one can construct a $2$-dimensional super lattice model with the partition function 
 $$Z_R^{lat}(\Lambda_{n,m})=\sum_{x\in\C: F^{ n}(x)=x} x^m$$
where $F:\C\to \C$ is a polynomial map\footnote{This seems to be a new type of integrability
 in lattice models, different from the usual Yang-Baxter ansatz.}, e.g. $F(x)=x^2+c$.  

The conclusion is that two different  proposals concerning motivic
local systems in positive characteristic:
 the first  (algebraic dynamics) and the third one
 (lattice models) are ultimately related. It is enough to find the
  dynamical realization, and then the lattice model will pop out.
  As it was mentioned already, most probably these two proposals would
  fail, but they still can serve as sources of analogies.

\vspace{2mm}

IHES, 35 route de Chartres, Bures-sur-Yvette 91440, France

{maxim@ihes.fr}

\end{document}